\documentclass{rspublic}

\usepackage{amsfonts}
\usepackage{amssymb}
\usepackage{graphicx}

\begin{document}

\title[Exponential asymptotics and boundary value problems]{Exponential asymptotics and boundary value problems: keeping both sides happy at all orders.}
\author[C.J.Howls]{C.J.Howls}
\affiliation{School of Mathematics, University of Southampton, Southampton SO17 1BJ, UK}
\label{firstpage}
\maketitle

\begin{abstract}{exponential asymptotics, asymptotics, hyperasymptotics, transseries, boundary value problems.}

We introduce templates for exponential asymptotic expansions that, in contrast to matched asymptotic approaches, enable the simultaneous satisfaction of both boundary values in classes of linear and nonlinear equations that are singularly perturbed with an asymptotic parameter $\epsilon\rightarrow 0^+$ and have a single boundary layer at one end of the interval.  For linear equations, the template is a transseries that takes the form of a sliding ladder of exponential scales.  For nonlinear equations, the transseries template is a two-dimensional array of exponential scales that tilts and realigns asymptotic balances as the interval is traversed.  An exponential asymptotic approach also reveals how boundary value problems force the surprising presence of transseries in the linear case and negative powers of $\epsilon$ terms in the series beyond all orders in the nonlinear case.  We also demonstrate how these transseries can be resummed to generate multiple-scales-type approximations that can generate uniformly better approximations to the exact solution out to larger values of the perturbation parameter.  Finally we show for a specific example how a reordering of the terms in the exponential asymptotics can lead to an acceleration of the accuracy of a truncated expansion.

\end{abstract}

\section{Introduction}

In this paper we examine the role of exponential asymptotics in boundary value problems (BVPs) for classes of second-order ordinary differential equations with a small positive asymptotic parameter $\epsilon$ multiplying the highest derivative.  In particular we consider both linear and nonlinear problems with boundary layers at one end of the finite domain.  At first sight, this appears to be a well-trodden classical problem with little more to be revealed.  Systems that cannot be solved exactly can be attacked using a vast toolbox of asymptotic approaches, most notably matched asymptotic expansions (MAE), see for example 
Nayfeh (1973) and references therein.  

Despite, or perhaps because of, the success of MAE approaches, BVPs have received less attention from asymptotic analysts than initial value problems.   Results on existence and uniqueness of solutions of second-order nonlinear two-point BVPs predated MAE 
(Wasow 1956 and references therein).  Formal theorems on the existence and character of asymptotic solutions of BVPs have been proved 
(Wasow 1956; Howes 1978). More recent asymptotic BVP work (e.g., O'Malley \& Ward 1997;  Ou \& Wong 2003, 2004; Ward 2005) has considered shock fronts and multiple spike solutions.   However, we shall show that, even for such a well-studied problem, the deployment of an exponential-asymptotic approach throws up subtle, surprising and generic asymptotic issues left unanswered by previous approaches.  

Exponential asymptotics (Berry 1989; Segur {\it et al} 1991) seeks to obtain a formal expansion of a function that incorporates all relevant exponential scales.  These scales are manifested typically by the exponential prefactors of (often) divergent infinite series.  Techniques exist for managing the divergence and analytical continuation of these series (Berry \& Howls 1991; Olde Daalhuis 1993), which lead to exponentially improved numerical approximations that are uniformly valid in wider ranges of the underlying parameters 
 (Olde Daalhuis 1998, Olde Daalhuis \& Olver 1998). 
This approach can provide better error bounds (Boyd 1994) and algebraic methods to resolve complex geometric structures (Howls 1997;  Delabaere \& Howls 2002).
In this paper we shall not be concerned with the implementation of the hyperasymptotics.  Rather, we seek to identify the templates, or form for the expansions that are required as input to these hyperasymptotic approaches.  The templates are the atomic building blocks which can be manipulated, using resummation techniques, to produce more sophisticated, `molecular asymptotic' approximations.  

Exponential asymptotic approaches to solutions of linear and nonlinear ordinary differential equations have been well explored 
(Chapman {\it et al} 1998; Olde Daalhuis 1998, 2004, 2005{\it a}{\it b}; Howls \& Olde Daalhuis 2003; \'{ O}lafsd\'{o}ttir {\it et al} 2005),
but have concentrated on initial value problems.   O'Malley and Ward (1997) 
applied exponential asymptotics to linear convection-diffusion PDE BVPs.  Lee \& Ward (1995)
considered exponentially sensitive ill-conditioning in nonlinear ODE BVPs.  Ward (2005)
also provided an asymptotic description of multiple spike solutions and the associated interior boundary layers in Carrier-Pearson-type systems. However this paper is the first time that a detailed analysis of the satisfaction of singularly perturbed BVPs at each and every exponential order has taken place.  The approach reveals the following surprising results.

First,  as is well known, MAE usually violate one of the boundary conditions.  This error may be negligible at small values of the asymptotic parameter $\epsilon$, but as $\epsilon$ grows the range of validity (both numerical and analytical) of the expansion may decrease. Here we show how, for a similar effort, an exponential asymptotic approach satisfies both boundary conditions and hence obtains an expansion that is valid in a wider range with at least a comparable accuracy.  

Second, an exponential approach also reveals a generic, yet previously unknown, intricate and interlocking ladder (linear) or array (nonlinear) structure within the expansions that is essential to the satisfaction of boundary conditions.  

Third, for initial value problems arising from linear ordinary differential systems, usually only a finite number of asymptotic series play a role, each being a leading order (in the exponential-asymptotic sense) expansion of one of the solutions.  The analysis here demonstrates explicitly how boundary values force the presence of  ``transseries" in linear systems here being infinite sums of exponentially-prefactored asymptotic series in the small parameter $\epsilon$.  Transseries occur naturally in nonlinear differential systems (Olde Daalhuis 2005{\it a}{\it b}; Costin 1998; Costin \& Costin 2001; Chapman {\it et al} 2007) due to the mixing of exponential scales by nonlinear terms, but their role in linear systems has not been highlighted explicitly before.  They are significant, because the resummation of transseries can generate non-local information about the singularity structure of the problem and can be used to extend the range of validity of solutions.

Fourth, for nonlinear BVPs, exponential asymptotics uncovers the presence of, at first sight, paradoxical {\it negative} powers of the small parameter $\epsilon$ in expansions.  Fortunately these are contained within series that are premultiplied by decaying exponentials in $1/\epsilon$ and so behave regularly as $\epsilon \rightarrow 0^+$.

Fifth, there is a choice in the order in which terms are resummed and this may be exploited, where possible, to allow for the derivation of numerically more efficient approximations to the solutions of nonlinear BVPs.

At the outset, we stress that we are not proposing that an exponential-asymptotic approach should replace that of MAE (or any other) in boundary layer calculations.  MAE remains a powerful and (where appropriately applied) numerically accurate approach.  We aim to highlight the underlying asymptotic structure of BVPs to facilitate extensions to a matched asymptotic approach when needed.

The outline of the paper is as follows.  In section 2, we study a pedgagocial linear BVP to highlight the issues from an exponential asymptotic viewpoint.  In section 3 we use general arguments to introduce the presence of transseries and the associated sliding ladder of exponential scales, demonstrating in section 4 how this will be a generic phenomenon in linear BVPs.  In section 5 we introduce nonlinear BVPs, again with a pedagogical example, and the transseries array structure that allows boundary value satisfaction at both ends of the interval.  The coefficients in the leading exponential orders are derived in sections 6 and 7 with a demonstration of the intricate realignment of balancing of terms, together with a demonstration of how it is possible to resum them to generate multiple-scales solutions.  Moving to higher order terms in section 8 we uncover a contradiction that forces the presence of negative powers of the small parameter in the transseries.   In section 9 we discuss how the reordering of terms within transseries can accelerate numerical agreement.  We end with a discussion and suggestions for future work.

\section{A linear example of boundary failure}

We take a pedagogical approach and first illustrate the issues we wish to address through a specific example.  The general case will be explained later.  We consider
\begin{eqnarray}\label{1}
\epsilon u''(x)+(2x+1)u'(x)+2u(x)=0, \qquad 0<x<1, \\
u(0)=\alpha, \ \  u(1)=\beta,  \qquad 0<\epsilon<<1,
\end{eqnarray}
where prime denotes differentiation with respect to $x$.  Conventionally the argument for a matched asymptotic approach runs along the following lines. Substitution of the standard ansatz
\begin{equation}\label{simple}
u(x)\sim\sum_{r=0}^{\infty}a_r(x)\epsilon^r
\end{equation}
into (\ref 1) generates first order linear differential equations for $a_r(x)$.  Hence $a_r(x)$ cannot in general satisfy the two boundary conditions simultaneously.  A boundary layer exists and there is an apparent need for a second expansion that obeys a rescaled equation, valid in the boundary layer and which satisfies the inner boundary condition.  The inner expansion is then matched to the outer expansion involving $a_r(x)$ that satisfies the other boundary condition.  The inner and outer expansions are matched to determine the remaining unknown constants.  

Rescaling (\ref{1}) with a change of variables $x=\epsilon X$ 
we may derive a composite expansion (at leading order)
\begin{equation}\label{leadm}
u_{\rm matched}(x)= {3\beta}/{(1+2x})+(\alpha-3\beta)\re^{-x/\epsilon}.
\end{equation}
Clearly we have $u_{\rm matched}(0)=\alpha$.  However $u(1)=\beta+{\mathcal O}\left(\re^{-1/\epsilon}\right)$ and so fails to satisfy the boundary condition at $x=1$, by at least the size of the inner expansion.  For small values of $\epsilon$, this may well be a numerically negligible error, but this error will grow with $\epsilon$, so reducing the range of validity of the solution.  


In this example, we learn more by following the approach of Latta (Nayfeh 1973, pp.145-154).  We try the simple exponential asymptotics approach of a WKB ansatz 
\begin{equation}\label{latta}
u(x)\sim \sum_{r=0}^{\infty}a_r(x)\epsilon^r+\re^{-F(x)/\epsilon}\sum_{r=0}^{\infty}b_r(x)\epsilon^r.
\end{equation}
We shall require $F(0)=0$ so that terms in the two series can balance at ${\mathcal O}(\epsilon^r)$ when $x=0$.  However, for general $x\ne 0, 1$, substitution of (\ref{latta}) into (\ref{1}) and balancing at each order ${\mathcal O}(\epsilon^r)$ and ${\mathcal O}\left(\re^{-F(x)/\epsilon}\right)$, we obtain:
\begin{eqnarray}
{\mathcal O}(\epsilon^0): &\qquad& (2x+1)a_0'(x)+2a_0(x)=0, \\
{\mathcal O}(\epsilon^1): &\qquad& (2x+1)a_1'(x)+2a_1(x)=-a_0''(x),  \\
{\mathcal O}\left(\re^{-F(x)/\epsilon}\epsilon^{-1}\right): &\qquad& F'(x)^2-(2x+1)F'(x)=0.
\end{eqnarray}
The last equation coupled with $F(0)=0$ gives  
\begin{equation}\label{Feqn}
F(x)=x^2+x.
\end{equation}

$F(x)$ is positive on $0<x<1$, as it should be, by assumption, so that (at least from an exponential-asymptotic viewpoint) no turning point occurs on that interval that can give rise to an interior boundary layer.  The remaining equations for $a_r$ and $b_r$ can be solved order by order.  The simplified equations for $b_r(x)$ are given below.   

The boundary conditions at $x=0$ are
\begin{equation}\label{bc0}
a_r(0)+b_r(0)=\delta_{r0}\alpha,
\end{equation}
where $\delta_{ij}$ is the Kronecker delta.  These can be satisfied by design.
However the boundary condition at $x=1$ can still only be satisfied up to exponential accuracy in $\epsilon$:
\begin{equation}\label{bc1}
a_r(1)+{\mathcal O}\left(\re^{-F(1)/\epsilon}\right)=\delta_{r0}\beta.
\end{equation}
Ignoring the exponential error at the right-hand boundary we generate the relations:
\begin{equation} \label{aterm}
a_r(x)=(a_{r-1}'(1)-a_{r-1}'(x))/(2x+1), \qquad a_0(x)=3\beta/(2x+1),
\end{equation}
\begin{equation}  \label{bterm}
b_r'(x)={b_{r-1}''(x)}/({2x+1}), \qquad b_r(0)=-a_r(0), \qquad b_0(x)=\alpha-3\beta.
\end{equation}
The latter equations show that the $b_r$ are here actually constants.  Hence the leading order solution is
\begin{equation}\label{approx}
u_{\rm WKB}(x)= \left({3\beta}/({2x+1})+(\alpha-3\beta)\re^{-(x^2+x)/\epsilon}\right)\left(1+{\mathcal O}\left(\epsilon\right) \right).
\end{equation}
Note that for $x>0$, the exponentially small term could be absorbed into the ${\mathcal O}(\epsilon)$ term, but is kept separate to delineate the existence of the boundary layer.

The matched and Latta-WKB approaches provide asymptotic ``solutions" that satisfy the boundary condition near to the boundary layer at $x=0$.   They both fail to satisfy the boundary condition at $x=1$. From an exponential asymptotic  (and for $x(x+1)/\epsilon={\mathcal O}(1)$, a numerical) viewpoint, this is unsatisfactory.  Comparison of  (\ref{leadm}) and (\ref{approx}) shows that, for a similar amount of effort, we have
$u_{\rm matched}(1)= \left(\beta+(\alpha-3\beta)\re^{-1/\epsilon}\right)\left(1+{\mathcal O}\left(\epsilon\right) \right),$
$u_{\rm WKB}(1)= \left(\beta+(\alpha-3\beta)\re^{-2/\epsilon}\right)\left(1+{\mathcal O}\left(\epsilon\right) \right).$
The order of the exponential error in the WKB approach at the right-hand boundary is here less than in the matched case.   This suggests that an extension of the WKB approach may be worth pursuing.  
We develop this approach in the next section.

\section{A transseries approach}

The exponential prefactors of the $a_r$-series and $b_r$-series in (\ref{latta}) balance at $x=0$ by design.  At $x=1$, the $a_r$ satisfy the boundary data, but the $b_r$ do not.  

Now suppose there is an additional series present, satisfying the differential equation (\ref{1}), but prefactored by $\re^{-F(1)/\epsilon}$.  At $x=1$, this can be used to cancel off the contribution from the $b_r$-series.  However, while addressing the satisfaction of the boundary condition at $x=1$, the presence of this new series means that the boundary condition at $x=0$ is now violated exponentially.  

This new series can still satisfy the boundary conditions at $x=0$ if there is another, additional, series prefactored by $\re^{-(F(x)+F(1))/\epsilon}$. For then, at $x=0$ with $F(0)=0$, this second new series balances the first new series at ${\mathcal O}(e^{-F(1)/\epsilon})$.  

\begin{figure}[htbp]
       \begin{center}
 \includegraphics[scale=0.425, angle=-90]{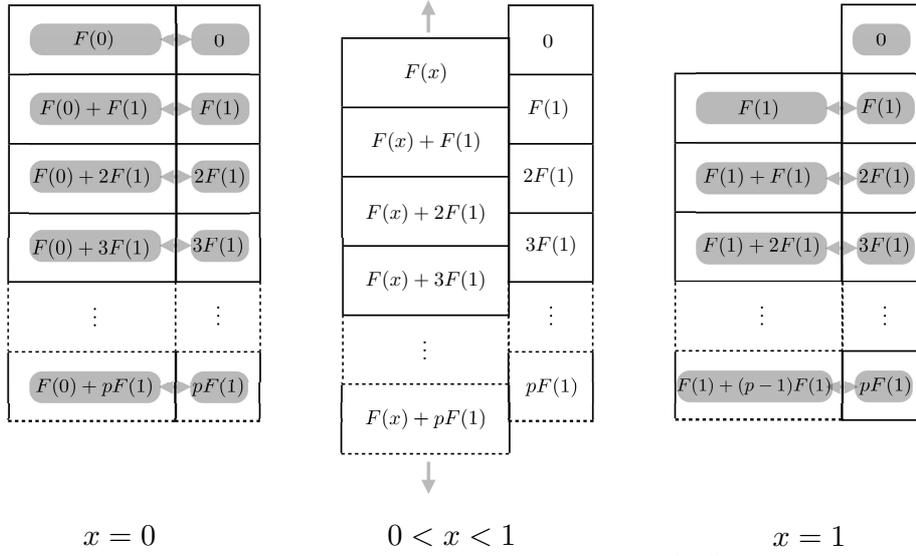}
   \caption{Ladder structure of exponentials within the transseries (\ref{template}).  The $0$ on the right of each panel denotes the series (\ref{simple}). The $pF(1), p=1,2,3 \dots$ on the right-hand side of each panel denotes the series prefactored by $\exp(-pF(1)/\epsilon)$.  For general $x$ the transseries exponents $F(x)+pF(1)$, $p=0, 1,2,3 \dots$ on the left translate up and down so that at each boundary there is a balance available at each exponential order. }
      \label{walk}
       \end{center}
\end{figure}

In turn, this second series can satisfy the boundary condition at $x=1$ if there is another, third additional series prefactored by $\re^{-2F(1)/\epsilon}$.  In turn this can satisfy the $x=0$ condition if there is another, fourth additional series prefactored by $\re^{-(F(x)+2F(1))/\epsilon}$  and so on and so forth.   (Problems with a single boundary layer near $x=1$ would require $F(1)$=0 in (\ref{Feqn}) and a reversal of roles of $F(0)$ and $F(1)$.) In summary, the boundary values generate a ladder of series with exponential prefactors defining the runs (see figure \ref{walk}). 

The middle column in figure \ref{walk} can be visualised as a sliding scale of prefactors that can be used
to line up different exponential scales.  At $x=0$ the scale lines up with both with the $0$-prefactor and also the scales at $pF(1)$, for integer $p\ge0$ (left hand panel). This ensures that the boundary conditions at $x=0$ are satisfied.  As $x$ travels from $0$ to $1$ this scale slides down (middle panel) so as to line up with the exponential prefactors $pF(1)$, $p\ge1$ at $x=1$ (right-hand panel).  The terms that are algebraic in $\epsilon$ in the re-balanced transseries can then be paired up so as to satisfy the boundary conditions, as we show below.  With this ladder, the zero-prefactored series exactly satisfies the boundary condition at $x=1$, with no exponential error.  

The template for the full asymptotic expansion is thus:
\begin{equation}\label{template}
u_{\rm trans}(x;\epsilon)\sim\sum_{p=0}^{\infty}\re^{-pF(1)/\epsilon}\sum_{r=0}^{\infty}a_r^{(p)}(x)\epsilon^r+\sum_{p=0}^{\infty}\re^{-\left(F(x)+pF(1)\right)/\epsilon}\sum_{r=0}^{\infty}b_r^{(p)}(x)\epsilon^r,
\end{equation}
where $a_r^{(0)}(x)=a_r(x)$.   Such a sum of exponentially-prefactored series is called a transseries 
(Olde Daalhuis 2005{\it a}, {\it b}; Costin 2001; Chapman {\it et al} 2007). The $p$-sums are transseries with prefactors that only depend on $\epsilon$, not $x$.   When approached directly from the asymptotics, the presence of transseries in linear equations is, at first sight, surprising.

We shall show below that this is the general form of the exponential-asymptotic template for the class of problems under consideration.  
More general templates might be derived for higher order equations.  The template for BVPs with internal boundary layers must be modified to take account of the associated connection problems.

Substituting (\ref{template}) into (\ref{1}) and hereafter using $\cdot_{,x}$ to denote $x$-derivatives, we arrive at the recurrence relations 
\begin{eqnarray}
a_r^{(p)}(x)&=&\left(a_{r-1,x}^{(p)}(1)-3b_r^{(p-1)}(1)-a_{r-1,x}^{(p)}(x)\right)/\left(2x+1\right), \\
b_{r,x}^{(p)}(x)&=&b_{r-1,xx}^{(p)}(x)/(2x+1).
\end{eqnarray}
with $b_r^{(-1)}(x)=0$.
The BVP then becomes
\begin{equation}\label{pbcs0}
a_r^{(p)}(0)+b_r^{(p)}(0)=\delta_{r0}\delta_{p0}\alpha, 
\qquad a_r^{(p)}(1)+b_r^{(p-1)}(1)=\delta_{r0}\delta_{p0}\beta.
\end{equation}
Simple calculations give the leading orders of each series as
\begin{equation}
a_0^{(0)}(x)=3\beta/(2x+1), \ \ \ a_0^{(p)}(x)=-3^p(\alpha-3\beta)/(2x+1), \ \ \ b_0^{(p)}(x)=3^p(\alpha-3\beta).
\end{equation}
Hence we have
\begin{eqnarray}
u_{\rm trans}(x;\epsilon)&=&\left\{\frac{3\beta}{2x+1}-\frac{(\alpha-3\beta)}{2x+1}\sum_{p=1}^{\infty}3^p \re^{-pF(1)/\epsilon} \right.  \nonumber \\  
&\ & \left.  + \  (\alpha-3\beta)\sum_{p=0}^{\infty}3^p \re^{-(F(x)+pF(1))/\epsilon}\right\}(1+{\mathcal O}(\epsilon)).
\end{eqnarray}
The $p$-sums are convergent for $|3\re^{-F(1)/\epsilon}|<1$.  Summing these is equivalent to changing the order of the $p$ and $r$ sums in transseries template (\ref{template}),  
(Olde Daalhuis 2005{\it a}, {\it b}; Costin 2001; Chapman {\it et al} 2007). The result is
\begin{equation}\label{psum}
u_{\rm trans}(x;\epsilon)=\left(\frac{3}{2x+1}\left\{\frac{\beta-\alpha \re^{-F(1)/\epsilon}}{1-3\re^{-F(1)/\epsilon}}\right\}+\left\{\frac{\alpha-3\beta}{1-3\re^{-F(1)/\epsilon}}\right\}\re^{-F(x)/\epsilon}\right)\left(1+{\mathcal O}(\epsilon)\right). 
\end{equation}

The approximation (\ref{psum}) satisfies both boundary conditions.  In this sense (at least) approximation (\ref{psum}) is better than either of the WKB or MAE approaches.  If terms in $\re^{-F(1)/\epsilon}$ are neglected in (\ref{psum}), we recover the WKB result (\ref{approx}).  

Higher order approximations can in principle be derived by similar resummations of the transseries at each order of ${\mathcal O}(\epsilon^r)$.  The results satisfy the boundary conditions at each order.  In Appendix A we show how the $p$-series may be resummed by a multiple scales approach.

\begin{figure}[htbp]
    \begin{center}
       \includegraphics[scale=0.425, angle=-90]{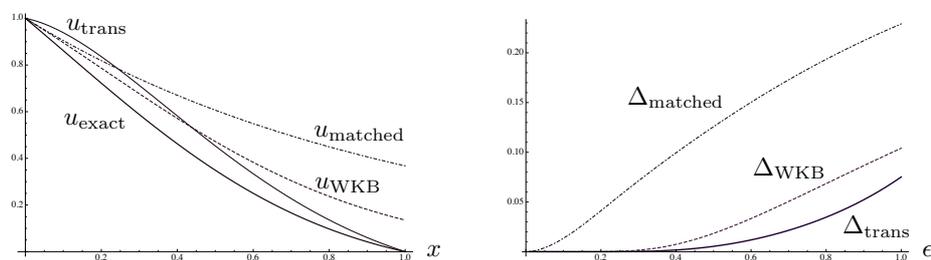}
   \caption{A comparison of the exact and various asymptotic solutions of (\ref{1}) for $\alpha=1$, $\beta=0$, $\epsilon=1$ over the range $0\le x\le 1$ (left).  The bold line is the exact solution, the thin solid line is the exponential asymptotic approximation, the dashed line is the WKB approximation, the dot-dashed line is the leading order MAE.   On the right the integrated leading order error (\ref{integerror}) over the range $0< x < 1$ as a function of $0<\epsilon\le 1$ is plotted.}
    \label{Asympplots}
    \end{center}
\end{figure}

In figure \ref{Asympplots} we compare the leading order results of the exponential asymptotic, Latta-WKB and matched asymptotic approaches with the exact solution for $\alpha=1$, $\beta=0$ and the comparatively large small parameter $\epsilon=1$.  Clearly the exponential asymptotic approach is the only one to satisfy both boundary conditions.  For the large value of $\epsilon=1$, the other asymptotic approaches (especially MAE) agree noticeably less well with the exact result away from the boundary layer near $x=1$.  

In the right-hand graph in figure \ref{Asympplots} we  compare the integral of the relative errors over the entire range $0\le x\le 1$:
\begin{equation}\label{integerror}
\Delta_{\rm approx}(\epsilon)=\int_0^1 |u_{\rm exact}(x;\epsilon)-u_{\rm approx}(x;\epsilon)|\rd x,
\end{equation}
for values of $\epsilon$ between 0 and 1.  Clearly the cumulated error in the exponential asymptotic approach is smaller over the entire range of $x$.


\section{Comparison with the exact solution and the general case}

A similar analysis can be carried out for more general solutions of second-order linear two-point boundary value problems.  Consider the system:
\begin{eqnarray}\label{gen}
\epsilon u''(x)+c(x)u'(x)+d(x)u(x)=0, \qquad 0<x<1, \\
u(0)=\alpha, \qquad u(1)=\beta, \qquad 0<\epsilon <<1,
\end{eqnarray}
where $c(x), d(x)>0$ for $0<x<1.$  Let the solutions of the equation be 
\begin{equation}\label{us}
u_1(x;\epsilon)=T_1^{(x)}(\epsilon), \ \ u_2(x;\epsilon)=\re^{-F(x)/\epsilon}T_2^{(x)}(\epsilon), 
\end{equation}
where the $T_j^{(\zeta)}(\epsilon)$ and $F(x)$ are given by 
\begin{equation}
T_j^{(\zeta)}(\epsilon)=\sum_{r=0}^{N-1}\tilde{T}_{j,r}^{(\zeta)}\epsilon^r+R_{j,N}(\zeta,\epsilon), \qquad F(x)=\int_0^xc(\eta)d\eta. 
\end{equation}
For a choice of $u_1$ and $u_2$ satisfying (\ref{gen}), the self-consistent solution of the BVP is
\begin{equation}\label{exactu}
u(x;\epsilon)=\frac{\left[\alpha u_2(1)-\beta u_2(0)\right]u_1(x)+\left[\beta u_1(0)-\alpha u_1(1)\right]u_2(x)}{\left[u_1(0)u_2(1)-u_1(1)u_2(0)\right]},
\end{equation}
provided the denominator does not vanish.
This can be written as
\begin{eqnarray}\label{gen1}
u(x;\epsilon)=&\ &\left\{\frac{\alpha-\beta T_1^{(0)}(\epsilon)/T_1^{(1)}(\epsilon)}{1-\re^{-F(1)/\epsilon}T_1^{(0)}(\epsilon)T_2^{(1)}(\epsilon)/T_1^{(1)}(\epsilon)T_2^{(0)}(\epsilon)}\right\}\frac{T_2^{(x)}(\epsilon)}{T_2^{(0)}(\epsilon)} \re^{-F(x)/\epsilon} \nonumber \\
&+& \left\{\frac{\beta-\alpha \re^{-F(1)/\epsilon}T^{(1)}_2(\epsilon)/T^{(0)}_2(\epsilon)}{1-\re^{-F(1)/\epsilon}T_1^{(0)}(\epsilon)T_2^{(1)}(\epsilon)/T_1^{(1)}(\epsilon)T_2^{(0)}(\epsilon)}\right\}\frac{T_1^{(x)}(\epsilon)}{T_1^{(1)}(\epsilon)}.
\end{eqnarray}
Replacing the $T_j^{(\zeta)}(\epsilon)$ by their asymptotic expansions and expanding the quotient leads naturally to transseries (\ref{template}).  Hence for linear second-order boundary value systems of type (\ref{gen}) the ladder transseries template (\ref{template}) is the general form.   

\section{Nonlinear second-order 2-point boundary value problems}

Given the ladder template for linear second-order 2-point BVP, it is natural to ask if a generalisation exists for analogous nonlinear boundary value systems.

Such systems have been studied asymptotically for decades and form the raison-d'\^{e}tre for many matched asymptotic approaches to boundary layers.  Significant results were provided by  Wasow (1956), O'Malley (1969), Howes (1978).
Wasow (1956)
established conditions under which systems of the form
\begin{equation}\label{Wasoweq}
\epsilon u''(x)=F_1(u,x)u'(x)+F_2(u,x), \qquad u(0)=\alpha, \ u(1)=\beta, \qquad 0<\epsilon<<1,
\end{equation}
have absolutely convergent perturbative series on the interval $0<x<1$.  
 O'Malley (1969)
considered a systems approach to the problem, and Howes (1978)
used the stability of the $\epsilon=0$ problem to study the existence of boundary, shock and corner layer solutions.  

We shall again restrict ourselves to a situation where only a single boundary layer exists, so follow Howes (1978) p.79 and consider his example (E3):
\begin{equation}\label{node}
\epsilon u''(x)+u'(x)u(x)-u(x)=0, \qquad u(0)=\alpha, \ u(1)=\beta,
\end{equation}
where $1<\beta<\alpha+1$. Under these conditions (Howes 1978, p.83) a boundary layer exists near to $x=0$.  Generalisations of what follows could include, for example,  larger regions of the $(\alpha, \beta)$ plane or $\epsilon$-dependent boundary data.  

The conventional MAE approach to (\ref{node}) generates the composite expansion 
\begin{equation}\label{unonmatched}
u_c(x)=\left(x+(\beta-1)\frac{(\beta-1)\tanh\left(\frac{(\beta-1)x}{2\epsilon}\right)+\alpha}{\beta-1+\alpha\tanh\left(\frac{(\beta-1)x}{2\epsilon}\right)}\right)\left(1+{\mathcal O}\left(\epsilon\right)\right).
\end{equation}
Note that while $u_c(0)=\alpha$, we again have $u_c(1)=\beta+{\mathcal O}\left(\re^{-(\beta-1)/\epsilon}\right)\ne \beta$.
Clearly, when $\beta>1$,  this error may diminish as $\epsilon \rightarrow 0^+$.  However for larger values of $\epsilon$, or small positive $\beta-1$ this error may become numerically significant.

The WKB approximation to the solution takes the form:
\begin{equation}\label{lattesoln}
u_{\rm WKB}(x)=\left((x+\beta-1)+\frac{(1+\alpha-\beta)(\beta-1)}{x+\beta-1}\re^{-x\left(\beta-1+x/2\right)/\epsilon}\right)\left(1+{\mathcal O}\left(\epsilon\right)\right).
\end{equation}
Clearly we have $u_{\rm WKB}(0)=\alpha$, but $u_{\rm WKB}(1)=\beta+{\mathcal O}\left(\re^{-\left(\beta-1/2\right)/\epsilon}\right)$.
Hence, as with MAE, the WKB template also fails to satisfy the boundary data at $x=1$.  



The WKB approach is just the first two terms of a transseries expansion of the general solution of nonlinear ODE which takes the form 
\begin{equation}
u(x)\sim \sum_{n=0}^{\infty}C^nu_{n}(x,\epsilon), \qquad u_n(x, \epsilon)\sim \re^{-nF(x)/\epsilon}\sum_{r=0}^{\infty}a_r^{(n)}\epsilon^r, \qquad F(0)=0,
\end{equation}
with the constant $C$ determined by the boundary conditions.  It will here be absorbed into the leading order coefficient of $u_1$, namely $a_0^{(1)}$.  The necessity to include integer powers of $e^{-F(x)/\epsilon}$ is caused by the nonlinear terms in the differential equation which generate successively higher order exponential scales that must be balanced (at each $x\ne 0$), at least in a formal solution.

However, we see immediately that this template still suffers from the same problem as the previous approaches.  Specifically the boundary conditions 
then require:
\begin{equation}\label{cap3}
\sum_{n=0}^{\infty}a_r^{(n)}(0)=\delta_{r0}\alpha, \qquad a_{r}^{(0)}(1)+\sum_{n=1}^{\infty}\re^{-nF(1)/\epsilon}a_r^{(n)}(1)``="\delta_{r0}\beta, 
\end{equation}
at each $r=0, 1, 2, \dots$.
For the time being, we put aside the issue of being able to satisfy the boundary conditions with an infinite sums of terms.  The point is that, for $F(1)\ne 0$, the boundary conditions at $x=1$ still leave terms unbalanced at ${\mathcal O}\left(\re^{-nF(1)/\epsilon}\right)$ for any integer value of  $n>0$. 

Inspired by the linear example above, we modify the transseries template as
\begin{equation}\label{cap4}
u(x)\sim\sum_{n=0}^{\infty} \sum_{p=0}^\infty u_n^{(p)}(x,\epsilon),
\end{equation}
\begin{equation}\label{cap55}
u_0(x)\sim \sum_{r=0}^{\infty}a_r^{(0,0)}(x)\epsilon^r, \qquad u_n^{(p)}(x,\epsilon)\sim \re^{-(nF(x)+pF(1))/\epsilon}\sum_{r=r_{\rm min}(n,p)}^{\infty}a_r^{(n,p)}(x)\epsilon^r. 
\end{equation}
The exponent $F(x)$ again slides between $F(0)\equiv0$ and $F(1)$ as $x$ runs between $0$ and $1$. The initial assumption of the minimum algebraic order in $\epsilon$, is $r_{\rm min}(n,p)=0$.  However, we shall see later that even in this simple example, this is not always so.


In contrast to the single exponent  $\re^{-F(x)/\epsilon}$ in the linear case, we now have an infinite number of exponential scales $\re^{-nF(x)/\epsilon}$, $n=0,1,2 \dots$.  At $x=0$ these all degenerate to $e^{-nF(0)/\epsilon}=1$ and so can be made to balance at each ${\mathcal O}(\epsilon^r)$ with the corresponding terms from $u_0$.  At $x=1$ these scales become  $\re^{-nF(1)/\epsilon}$, $n=1,2 \dots$.  To satisfy the boundary condition at each exponential order thus again requires the inclusion of additional series prefactored with scales $\re^{-pF(1)/\epsilon}$, $p=1,2 \dots$.  However, substitution shows that these boundary series will multiply series in $\re^{-nF(x)/\epsilon}$ in the $uu'(x)$ term to produce further contributions with prefactors $\re^{-(nF(x)+pF(1))/\epsilon}$, $n, p=0,1,2 \dots$.  All of these series have to be balanced at $x=0,1$ to satisfy the boundary data.  Fortunately this is possible, as explained in the diagrams below.

\begin{figure}[htbp]
       \begin{center}
\includegraphics[scale=0.45, angle=-90]{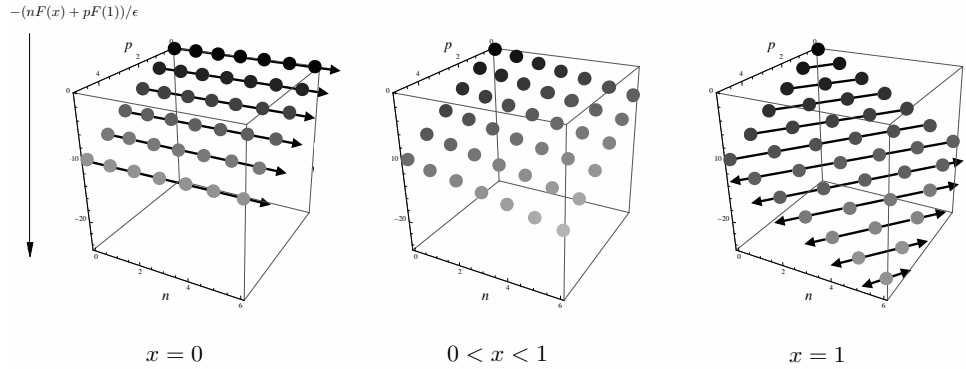}
    \caption{Array structure of exponentials within the double transseries (\ref{cap4}-\ref{cap55}). Dots joined by lines denote exponential scales $\re^{-(nF(x)+pF(1))/\epsilon}$ that balance at the value of $x$ shown.    Arrows denote the presence of additional dots corresponding to additional balancing exponential scales outside the bounding box of the diagrams.}
    \label{nonlineararray}
       \end{center}
\end{figure}

From (\ref{cap55}) the  exponents labelled by integers $(n,p)$ generate not a ladder of exponential scales, but a planar lattice (see figure \ref{nonlineararray}).  The heights of the lattice points are at the exponential scales $\re^{-(nF(x)+pF(1))/\epsilon}$.  
At $x=0$, $nF(0)=0$, and so balances occur between $n=0,1,2,3\dots$ series prefactored by $\re^{-pF(1)/\epsilon}$ at
 each value of $p$.   As $x$ increases to 1, $F(x)$ varies and the planar lattice tilts.  When 
 $x=1$, the prefactors become  $\re^{-(n+p)F(1)/\epsilon}$ and the scales have reconfigured so 
 that balancing now takes place between series prefactored by exponentials with $n+p=K$ for fixed $K$ 
 (along diagonal rows on the right in figure 3). We shall demonstrate this reconfiguration in detail by explicit 
 calculation for the system (\ref{node}) in the following sections.

\section{Derviation of the coefficients}

Substitution of (\ref{cap4}) into (\ref{node}) and balancing at ${\mathcal O}\left(\re^{-(nF(x)+pF(1))/\epsilon}\right)$ gives
\begin{equation}\label{cap6}
\epsilon u_{n,xx}^{(p)}(x)+\sum_{m=0}^{n}\sum_{q=0}^p u_m^{(q)}(x)u_{n-m,x}^{(p-q)}(x)-u_n^{(p)}(x)=0.
\end{equation}
For fixed $(n,p)$ we can substitute (\ref{cap55}) into  (\ref{cap6}) and balance at ${\mathcal O}(\epsilon^r)$.  We obtain the following balances at $(n,p)=(0,0)$, ${\mathcal O}(\epsilon^r)$ and for $n\ne 0$, ${\mathcal O}(1/\epsilon)$, respectively
\begin{equation}\label{cap7aab}
 \qquad a_{r-1,xx}^{(0,0)}(x)+\sum_{s=0}^r a_s^{(0,0)}(x)a_{r-s,x}^{(0,0)}(x)-a_r^{(0,0)}(x)=0.
\end{equation}
\begin{equation}\label{cap7}
n^2F'(x)a_0^{(n,p)}(x)-\sum_{m=0}^{n-1}\sum_{s=0}^p (n-m)a_0^{(m,s)}(x)a_0^{(n-m,p-s)}(x)=0.
\end{equation}

These recurrence relations now generate the leading order asymptotic terms, with associated arbitrary constants $C_r^{(n,p)}$.  From $(n,p)=(0,0)$ at ${\mathcal O}(\epsilon^0)$ we find:
\begin{equation}\label{cap10}
a_0^{(0,0)}(x)=x+C_0^{(0,0)}.
\end{equation}
At ${\mathcal O}(\epsilon)$ when $(n,p)=(0,0)$, we have
\begin{equation}\label{cap11}
a_{0,xx}^{(0,0)}+a_0^{(0,0)}a_{1,x}^{(0,0)}+a_1^{(0,0)}a_{0,x}^{(0,0)}-a_1^{(0,0)}(x)=0 \qquad \Rightarrow \qquad a_1^{(0,0)}(x)=C_1^{(0,0)}.
\end{equation}
It is then simple to prove from $(n,p)=(0,0)$, ${\mathcal O}(\epsilon^r)$, that 
\begin{equation}\label{cap12}
\qquad a_r^{(0,0)}(x)=C_r^{(0,0)}, \ r>0.
\end{equation}
$F(x)$ is deduced from $(n,p)=(1,0)$ at ${\mathcal O}(1/\epsilon)$:
\begin{equation}\label{cap8}
F'(x)a_0^{(1,0)}(x)-a_0^{(0,0)}(x)a_0^{(1,0)}(x)=0 \qquad \Rightarrow \qquad F'(x)=a_0^{(0,0)}(x).
\end{equation}
This assumes that $a_0^{(1,0)}(x)\ne0$.  This result allows us to recast (\ref{cap7}) as a recurrence relation to generate $a_0^{(n,p)}(x)$. 
For nonzero values of $n$ at ${\mathcal O}(1/\epsilon)$ we have:
\begin{eqnarray}\label{cap7a}
n(n-1)a_0^{(0,0)}a_0^{(n,p)}(x)&=&\sum_{m=0}^{n-1}\sum_{s=1}^p (n-m)a_0^{(m,s)}(x)a_0^{(n-m,p-s)}(x) \nonumber \\ &\ & + \sum_{m=1}^{n-1}(n-m)a_0^{(m,0)}a_0^{(n-m,p)}.
\end{eqnarray}
It is straightforward to show from (\ref{cap7a}) that 
\begin{equation}\label{cap7aa}
a_0^{(n,0)}=2a_0^{(0,0)}\left(a_0^{(1,0)}/2a_0^{(0,0)}\right)^n.
\end{equation}
The term $a_0^{(1,0)}(x)$ can be found, from $(n,p)=(1,0)$, ${\mathcal O}(\epsilon^0)$ and using (\ref{cap12}):
\begin{eqnarray}\label{cap7aaa}
-a_0^{(1,0)}(x)\left(1+a_0^{(0,0)}(x)C_1^{(0,0)}\right)-a_0^{(0,0)}(x)a_{0,x}^{(1,0)}&=&0,
\end{eqnarray}
\begin{equation}\label{cap7aaaa}
\Rightarrow \ a_0^{(1,0)}(x)=C_0^{(1,0)}\re^{-C_1^{(0,0)}(x+C_0^{(0,0)})}/(x+C_0^{(0,0)}).
\end{equation}

When $n=0$, $p\ne0$, further terms can be deduced to be constants from the ${\mathcal O}(\epsilon^0)$ equations by first observing that, when $p=1$:
\begin{equation}\label{cap20}
a_0^{(0,0)}(x)a_{0,x}^{(0,1)}+a_0^{(0,1)}(x)a_{0,x}^{(0,0)}(x)-a_0^{(0,1)}(x)=0, \qquad \Rightarrow \qquad a_0^{(0,1)}(x)=C_0^{(0,1)}.
\end{equation}
Substitution of this result into the corresponding equation for $p>1$ gives:
\begin{equation}\label{cap19}
\sum_{s=0}^p a_0^{(0,s)}(x)a_{0,x}^{(0,p-s)}(x)-a_0^{(0,p)}(x)=0.
\end{equation}
It follows by iteration that:  
\begin{equation}\label{cap21}
a_0^{(0,p)}(x)=C_0^{(0,p)}, \qquad p=1,2,3, \dots.
\end{equation}
We find these constants in the next section.


\section{Boundary conditions reconfiguring the exponential scales}
\begin{figure}[htbp]
       \begin{center}
        \includegraphics[scale=0.425, angle=90]{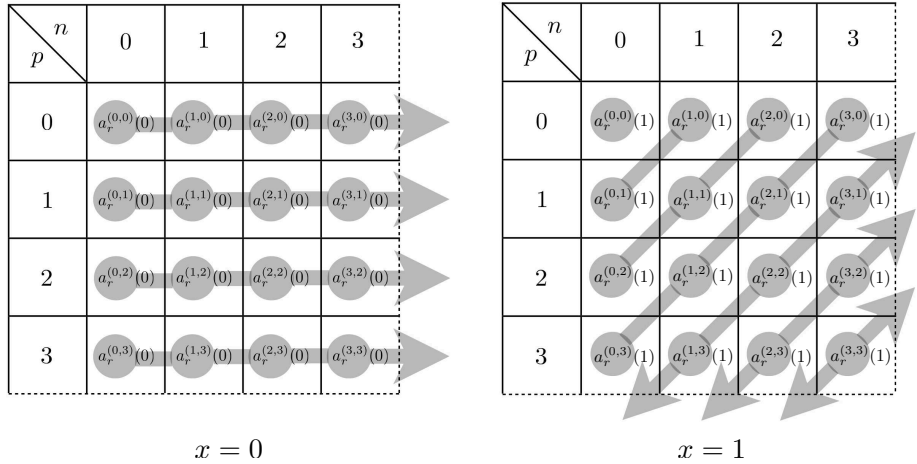}
    \caption{Balancing of terms $a_r^{(n,p)}$ at fixed $r$ between transseries with varying $(n,p)$ (\ref{cap4}-\ref{cap55}) at $x=0$ and $x=1$ due to the realignment of the exponential scales. Terms joined by the grey lines balance at the values of $x$ shown.  At $x=0$ at each fixed $p$, an infinite sum over all terms $n$ must satisfy the boundary conditions at the left hand end of the interval.  At $x=1$ terms with $n+p=K$ for each fixed $K=0,1,2 \dots$ must be summed to satisfy the boundary conditions at the right-hand end of the interval.}
    \label{termboundaryalign}
       \end{center}
\end{figure}

We now apply the boundary conditions, 
$u(0)=\alpha$, $u(1)=\beta,$
balanced first at each exponential order, and then at each algebraic order of $\epsilon$. The balancing of terms within the relevant exponential scales is illustrated graphically in figure \ref{termboundaryalign}.

In what follows, we set ${a_r^{(n,p)}}(0)={C_r^{(n,p)}}$.

At $x=0$, since $F(0)=0$ the exponential scales ${\mathcal O}\left(\re^{-nF(x)/\epsilon}\right)$ all collapse to the same degenerate level and the distinct exponential scales are given by ${\mathcal O}\left(\re^{-pF(1)/\epsilon}\right)$.  So for fixed $p$, the terms $a_r^{(n,p)}$ balance at order ${\mathcal O}\left(\re^{-pF(1)/\epsilon}\epsilon^r\right)$ according to:
\begin{equation}\label{capbc1}
\sum_{n=0}^{\infty}{a_r^{(n,p)}}(0)=\delta_{r0}\delta_{p0}\alpha, \ \ \Rightarrow \ \ \sum_{n=0}^{\infty}{a_0^{(n,0)}}(0)=\alpha, \ \ \sum_{n=0}^{\infty}{a_r^{(n,p)}}(0)=0, \ (p, r\ne0).
\end{equation}

At $x=1$, the template (\ref{cap4}-\ref{cap55}) reconfigures the exponential scales.  The orders ${\mathcal O}\left(\re^{-(n+p)F(1)/\epsilon}\right)$  with $n+p=K$, $K={\rm const}$, become degenerate.  Within these scales the terms $a_r^{(n,p)}$, consequently align at each algebraic ${\mathcal O}(\epsilon^r)$ according to:
\begin{equation}\label{capbc2}
\sum_{n=0}^{K}{a_r^{(n,K-n)}}(1)=\delta_{r0}\delta_{K0}\beta \ \ \Rightarrow \ \ {a_r^{(0,0)}}(1)=\delta_{r0}\beta, \ \ \sum_{n=0}^{K}{a_r^{(n,K-n)}}(1)=0, \ (n,r\ne0).
\end{equation}

We make progress by first considering the conditions at $x=1$.  Using (\ref{cap10}), (\ref{cap7aaaa}), (\ref{capbc2}) we can identify immediately that
\begin{equation}\label{capbc3}
a_0^{(0,0)}(1)= 1+C_0^{(0,0)} =\beta \qquad \Rightarrow \qquad a_0^{(0,0)}(x)= x+\beta-1,
\end{equation}
\begin{equation}\label{capbc3a}
a_r^{(0,0)}(1)= \delta_{r0}\beta=C_r^{(0,0)} =0 \qquad r>0, \ \ \Rightarrow \ \ a_0^{(1,0)}(x)=C_0^{(1,0)}/(x+\beta-1).
\end{equation}
Hence, returning to $x=0$ we have, for example,
\begin{equation}\label{capbc4}
\sum_{n=1}^{\infty}{a_0^{(n,0)}}(0)=\alpha-\beta+1.
\end{equation}
Using the form (\ref{cap7aa}) we may sum this infinite series formally as a geometric progression and obtain
\begin{equation}\label{capbc4a}
\frac{2a_0^{(0,0)}(0)a_0^{(1,0)}(0)}{2a_0^{(0,0)}(0)-a_0^{(1,0)}(0)}=\alpha-\beta+1,
\end{equation}
\begin{equation}\label{capbc4b}
\frac{2(\beta-1)C_0^{(1,0)}(0)}{2(\beta-1)^2-C_0^{(1,0)}(0)}=\alpha-\beta+1 \qquad \Rightarrow \qquad C_0^{(1,0)}=\frac{2(\alpha-\beta+1)(\beta-1)^2}{(\alpha+\beta-1)}.
\end{equation}

Returning back to $x=1$ we then have:
\begin{equation}
a_0^{(0,1)}(1)+a_0^{(1,0)}(1)=0 \ \ \Rightarrow \ \ C_0^{(0,1)}+\frac{C_0^{(1,0)}}{\beta}=0  \ \ \Rightarrow \ \ C_0^{(0,1)}=-\frac{2(\alpha-\beta+1)(\beta-1)^2}{\beta(\alpha+\beta-1)}.
\end{equation}
Hence we finally have the leading orders of the first three $(n,p)$ transseries as:
\begin{equation}\label{lead}
a_0^{(0,0)}(x)= x+\beta-1, \ \  a_0^{(1,0)}(x)= \frac{2(\alpha-\beta+1)(\beta-1)^2}{(\alpha+\beta-1)(x+\beta -1)},\ \ a_0^{(0,1)}(x)=-a_0^{(1,0)}(1).
\end{equation}
From (\ref{cap21}) and (\ref{capbc2}) we can also deduce that for this example $a_r^{(0,0)}(x)=0$, $r\ge1$.  Note that this means that for the special case of (\ref{node}), the expansion $u_0^{(0,0)}$ will not generate a Stokes phenomenon, since it truncates after the first term. This also reflects the fact that there are no internal boundary layers.  Note also this means that the exponent $F(x)$ is not here directly connected to a factorial-over-power late-term ansatz (Dingle 1973; Berry 1989; Chapman {\it et al} 1998).

We may resum the $n$-transseries for $p=0$ using (\ref{lead}) and obtain the following:
\begin{eqnarray}\label{nresum}
u_{\rm trans}(x)&=&\left(a_0^{(0,0)}(x)+\sum_{n=1}^{\infty}\re^{-nF(x)/\epsilon}a_0^{(n,0)}(x)\right)\left(1+{\mathcal O}(\epsilon)\right) \\ \label{nresum1}
&=&\left(a_0^{(0,0)}(x)+\frac{2a_0^{(0,0)}(x)a_0^{(1,0)}(x)\re^{-F(x)/\epsilon}}{2a_0^{(0,0)}(x)-a_0^{(1,0)}(x)\re^{-F(x)/\epsilon}}\right)\left(1+{\mathcal O}(\epsilon)\right).
\end{eqnarray}
For most practical purposes, this result may be all that is needed. 

\begin{figure}[htbp]
       \begin{center}
\includegraphics[scale=0.425, angle=-90]{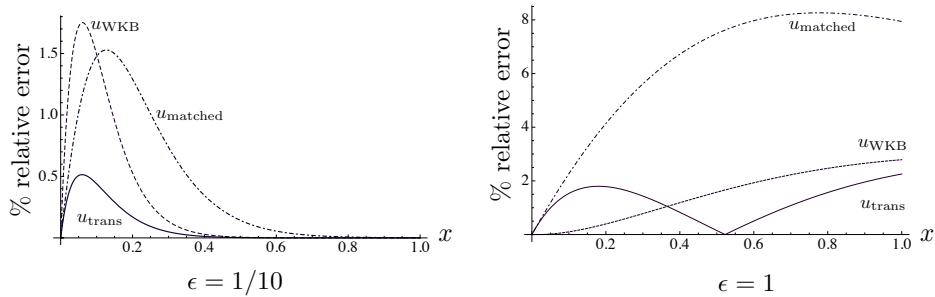}
    \caption{Comparison of the percentage relative errors, $100|1-u_{\rm approx}/u_{\rm exact}|$, for the matched (dotted), WKB (dashed), and resummed $n$-exponential transseries with no $p$-series (\ref{nresum}), (solid) for $\alpha=3/2, \beta=2$, with values of $\epsilon=1/10, 1$.}
    \label{comp1}
       \end{center}
\end{figure}

A comparison of the numerical errors for the MAE (\ref {unonmatched}), WKB (\ref{lattesoln}) and resummation (\ref {nresum}) is shown in figure \ref{comp1} for typical values of $\epsilon=1/10$ and $1$ with fixed $\alpha=3/2$ and $\beta=2$.  First, for small $\epsilon$ (left hand graph), although the errors committed by the approximations at $x=1$ are too small to display on the scale of the graph,  the transseries approximation is uniformly numerically better than either the matched or WKB solutions.  Secondly, for larger values of $\epsilon$ both the WKB and  transseries approximations are uniformly better approximations than the MAE.  The latter violates the boundary condition at $x=1$ by a noticeably larger amount.  Although it satisfies the boundary condition at $x=0$, the transseries approximation is not uniformly better than the WKB near there, agreeing more with the MAE.    However the transseries approximation quickly does beat WKB as $x$ increases and is definitely better at $x=1$.   
Note that the transseries approximation (\ref {nresum}) only contains terms that have $p=0$. If we summed the $n$-series with $p=1$ and added these in at ${\mathcal O}(\epsilon)$, the error at $x=1$ would diminish.  In the next section we show how the higher order balances reveal unexpected behaviour in the nonlinear template.

\section{Higher order balancing}

Moving on to higher exponential orders, at $x=1$ we need to balance terms involving the $(n,p)=(2,0), (1,1), (0,2)$ series (cf. figure \ref{termboundaryalign}):
\begin{equation}
a_0^{(0,2)}(1)+a_0^{(1,1)}(1)+a_0^{(2,0)}(1)=0.
\end{equation}
The behaviour of the $(n,p)=(2,0)$, $(0,2)$ terms have already been derived above:
\begin{equation}
a_0^{(2,0)}(x)=\frac{\left(a_0^{(1,0)}(x)\right)^2}{2a_0^{(0,0)}(x)}=\frac{2 (\alpha -\beta +1)^2(\beta -1)^4 }{(x+\beta -1)^3 (\alpha +\beta -1)^2}, \qquad a_0^{(0,2)}(x)=C_0^{(0,2)}.
\end{equation}
Seeking the $a_0^{(1,1)}(x)$ term we use (\ref{cap6}) with $(n,p)=(1,1)$:
\begin{equation}\label{n1p1}
\epsilon u_{1,xx}^{(1)}+u_0^{(0)}u_{1,x}^{(1)} +u_0^{(1)}u_{1,x}^{(0)} +u_1^{(0)}u_{0,x}^{(1)} +u_1^{(1)} u_{0,x}^{(0)}-u_1^{(1)}=0.
\end{equation}
Assuming $r_{\rm min}(1,1)=0$ (cf. (\ref{cap55})), at ${\mathcal O}(1/\epsilon)$ in (\ref{n1p1}), after simplification we have 
\begin{equation}
F'(x)a_0^{(0,1)}(x)a_0^{(1,0)}(x)=0.
\end{equation}
This obviously leads to a contradiction, since prior calculations show that neither of the conditions $F'(x)=0$, 
$a_0^{(0,1)}(x)=0$, or $a_0^{(1,0)}(x)=0$ hold.  

We may resolve this impasse by taking $r_{\rm min}(1,1)=-1$. For then, at ${\mathcal O}(1/\epsilon^2)$, (\ref{n1p1}) is satisfied by $F'(x)=a_0^{(0,0)}$ and at ${\mathcal O}(1/\epsilon)$ we have
\begin{equation}
a_0^{(0,0)}(x)a_{-1,x}^{(1,1)}(x)+a_{-1}^{(1,1)}(x)=-a_0^{(0,0)}(x)a_0^{(0,1)}(x)a_0^{(1,0)}(x),
\end{equation}
with general solution,
\begin{equation}\label{aminus1}
a_{-1}^{(1,1)}(x)=\frac{C_{-1}^{(1,1)} \beta  (\alpha +\beta -1)^2+4 x (\beta -1)^4 (\alpha -\beta +1)^2}{\beta  (x+\beta -1)
   (\alpha +\beta -1)^2}.
   \end{equation}
Given that, from above, there is no term $a_{-1}^{(0,1)}(0)$, the boundary condition that this term must satisify at $x=0$ is (cf. figure \ref{termboundaryalign}),
\begin{equation}\label{minus1bc}
\sum_{n=1}^\infty a_{-1}^{(n,1)}(0)=0.
\end{equation}
By examination of the ${\mathcal O}(1/\epsilon^2)$ balances in the $(n,p)=(1,p)$ equations (\ref{cap6}), it is possible to show
that 
\begin{equation}
a_{-1}^{(n,1)}(x)=n\left({a_0^{(1,0)}(x)/ 2 a_0^{(0,0)}(x)}\right)^{n-1}a_{-1}^{(1,1)}(x).
\end{equation}
Hence we can sum (\ref{minus1bc}) to obtain (cf. (\ref{aminus1}))
\begin{equation}
\frac{4 a_{-1}^{(1,1)}(0)\left(a_0^{(0,0)}(0)\right)^{2}}{\left(2a_0^{(0,0)}(0)-a_0^{(1,0)}(0)\right)^{2}}=0 \qquad \Rightarrow \qquad  a_{-1}^{(1,1)}(x)=\frac{4 x (\beta -1)^4 (\alpha -\beta +1)^2}{\beta  (x+\beta -1)
   (\alpha +\beta -1)^2}.
\end{equation}
Consequently, since $a_{-1}^{(2,0)}(x)=0$ we can now determine the constant $a_{-1}^{(0,2)}\equiv C_{-1}^{(0,2)}$ from the boundary condition at $x=1$:
\begin{eqnarray}
&a_{-1}^{(2,0)}(1)&+a_{-1}^{(1,1)}(1)+a_{-1}^{(0,2)}(1)=0  \nonumber \\  \Rightarrow \   &a_{-1}^{(0,2)}(x)&=- a_{-1}^{(1,1)}(1)=-\frac{4 (\beta -1)^4 (\alpha -\beta +1)^2}{\beta^2}.
\end{eqnarray}

Finally, it is then possible to resum the contributions from $a_{-1}^{(n,1)}$ to obtain
\begin{equation}\label{am1}
\sum_{n=0}^\infty \frac{a_{-1}^{(n,1)}(x)}{\epsilon}\re^{-(nF(x)+F(1))/\epsilon}=\frac{4 a_{-1}^{(1,1)}(x)\left(a_0^{(1,0)}(x)\right)^{2}\re^{-(F(x)+F(1))/\epsilon}}{\epsilon\left(2a_0^{(0,0)}(x)-a_0^{(1,0)}(x)\re^{-F(x)/\epsilon}\right)^{2}}.
\end{equation}
The exponential asymptotic approach has revealed the unexpected and counterintuitive presence of $1/\epsilon$ terms in the expansion as $\epsilon\rightarrow 0^+$.  However, these $1/\epsilon$ terms are prefactored by $\re^{-F(1)/\epsilon}$ with $F(1)>0$.  Hence they are not only beyond all algebraic orders, but also do vanish exponentially fast as $\epsilon\rightarrow 0^+$.   

The existence of an ${\mathcal O}(1/\epsilon)$ term in the $(n,p)=(1,1)$ series forces the presence of ${\mathcal O}(1/\epsilon^p)$ terms in series with $p>1$.  The corrected form of the template (\ref{cap55}) has:
\begin{equation}
r_{\rm min}(n,p)=\left\{\begin{array}{cc}-p, \qquad n>p, \\-{\rm Floor}[(n+p)/2], & n\le p.\end{array}\right.
\end{equation}
As above, since the negative powers of $\epsilon$ are multiplied by $\re^{-pF(1)/\epsilon}$, they will remain technically beyond all algebraic orders and vanish as $\epsilon\rightarrow 0^+$.  

The presence of these terms means that in generic cases, as $p$ increases, additional work must be done to obtain resummed expansions accurate to \\ ${\mathcal O}(\re^{-n(F(x)+F(1))/\epsilon}\epsilon)$.  For example, it is a simple algebraic exercise to show that the calculation to derive the resummation analogous to (\ref{am1}) of the ${\mathcal O}(\epsilon^0)$ terms with $p=1$, 
\begin{equation}\label{am2}
\sum_{n=0}^\infty \re^{-n(F(x)+F(1))/\epsilon}a_{0}^{(n,1)}(x)\epsilon^0,
\end{equation}
actually requires first a complete derivation of all the terms in $a_1^{(n,0)}(x)$.  

A detailed set of calculations generates the following results for $n\geq 1$:
\begin{equation}\label{Keqn}
a_{1}^{(n,0)}(x)=n\left({a_0^{(1,0)}/ 2 a_0^{(0,0)}}\right)^{n-1}a_{1}^{(1,0)}+{K_n/a_0^{(0,0)}}
\left({a_0^{(1,0)}/ 2 a_0^{(0,0)}}\right)^{n},
\end{equation}
\begin{eqnarray}\label{Leqn}
a_{0}^{(n,1)}(x)&=&n\left({a_0^{(1,0)}/ 2 a_0^{(0,0)}}\right)^{n-1}a_{0}^{(1,1)}   \nonumber \\ 
&\ &+ n(n-1){\left(a_0^{(1,0)}\right)^{n-2}/ \left(2a_0^{(0,0)}\right)^{n-1}}a_{1}^{(1,0)}a_{-1}^{(1,1)} \nonumber \\
&\ & -2(n-1)\left({a_0^{(1,0)}/ 2 a_0^{(0,0)}}\right)^{n}a_{0}^{(0,1)} \nonumber \\
& \ & +2L_n{\left(a_0^{(1,0)}\right)^{n-1}/ \left(2a_0^{(0,0)}\right)^{n+1}}a_{-1}^{(1,1)}.
\end{eqnarray}
The coefficients $K_n$ in (\ref{Keqn}) and $L_n$ in (\ref{Leqn}) satisfy the recurrence relations
\begin{eqnarray}
{(n-1)}K_n/2&=&\sum_{m=2}^{n-1}K_m-{(2n+1)(n-1)/ n},\\
{(n-1)}L_n/2&=&\sum_{m=2}^{n-1}L_m+\sum_{m=2}^{n-1}(n-m)K_m-(2n+1)(n-1).
\end{eqnarray}
A $z$-transform of these relations generates
\begin{equation}\label{Ksum}
\sum_{n=2}^\infty K_nt^n=H(t)\equiv{2t(1-2{\rm Li}_2(t))\over(1-t)^2}+{2(1+t)\ln(1-t)\over 1-t},
\end{equation}
\begin{equation}\label{Lsum}
\sum_{n=2}^\infty L_nt^n=G(t)\equiv{2t(1+t)(t-2{\rm Li}_2(t))\over(1-t)^3}+{8t\ln(1-t)\over (1-t)^2},
\end{equation}
where ${\rm Li}_2(t)$ is the dilogarithm function (DLMF 2010).

The sum over $n$ of (\ref{Keqn}) at $x=0$ can be achieved using (\ref{Ksum}) and according to the boundary value, be set equal to zero.  It only depends on unknown $a_1^{(1,0)}$ and the already known $a_0^{(1,0)}$ and $a_0^{(0,0)}$. Coupled with the general solution of 
\begin{equation}
a_0^{(0,0)}a_{1,x}^{(1,0)}+a_{1}^{(1,0)}=a_{0,xx}^{(1,0)},
\end{equation}
we can use this sum to determine $a_1^{(1,0)}(x)$.  This can be substituted into (\ref{Leqn}), which, in turn can be summed at $x=0$ using (\ref{Lsum}) to find $a_0^{(1,1)}(x)$.  In turn  (\ref{Leqn}) can then be multiplied by $\re^{-(nF(x)+F(1))/\epsilon}$ and summed over $n$ to obtain (\ref{am2}).  The resulting expression is large and so this is left as an exercise for the reader.

In principle this scheme can be repeated for each $p$.  For given $p$, the corresponding ${\mathcal O}(\epsilon)$ estimate will involve a sum over $n$ of $a_{0}^{(n,p)}(x)\re^{-(nF(x)+pF(1))}$.  To determine the $a_{0}^{(n,p)}(x)$ will likely require the determination of all terms from $a_{r}^{(n,p)}(x)$, $0>r>r_{\min}(n,p)$ beforehand.  This will become an increasingly arduous task.

However, as the transseries template contains three sums (in $n$, $p$ and $r$) there is the potential to alter the order of resummation of the terms to improve the numerical accuracy of the approximation.  

\section{Alternative resummations}

The above resummations are over $n$ for each $p$ and fixed order $\epsilon^r$.  Alternatively it may be possible to reorder the resummations, for example resumming first over $p$ and obtain numerically better agreement with the exact solution.

Due to the particular form of (\ref{node}) we may proceed here down this route by reverting to a single transseries formulation of the form
\begin{eqnarray}\label{pw1}
u(x)=\sum_{k=0}^{\infty}\lambda^k \tilde{u}_k(x; \epsilon), \qquad \tilde{u}_0(x, \epsilon)=\sum_{r=0}^{\infty} \tilde{a}_r^{(0)}(x)\epsilon^r,
\end{eqnarray}
where $\lambda$ is here just an ordering parameter that will be eventually set equal to one.  In the absence of the boundary data, the $\tilde{u}_k(x; \epsilon)$ can be regarded as providing terms of order ${\mathcal O}\left(\exp(-kF(x)/\epsilon)\right)$.
Substitution into (\ref{node}) and balancing at order $\lambda^k$ generates at ${\mathcal O}(\lambda^0)$:
\begin{eqnarray}\label{pw2}
\epsilon \tilde{u}_0''(x)+\tilde{u}_0(x) \tilde{u}_0'(x)-\tilde{u}_0(x)=0
\end{eqnarray}
and balancing at ${\mathcal O}(\epsilon^r)$ and setting $\tilde{a}_r^{(0)}(1)=\delta_{r0}\beta$ we have
\begin{eqnarray}\label{pw4}
\tilde{a}_r^{(0)}(x)=\delta_{r0}(x+\beta-1) \qquad \tilde{u}_0(x)=x+\beta-1,
\end{eqnarray}
as above.  Note that the truncated nature of $ \tilde{u}_0(x)$ is due to the particular system under study.  In general it would be only expressible as a formal, divergent series.

\begin{figure}[htbp]
       \begin{center}
        \includegraphics[scale=0.45, angle=90]{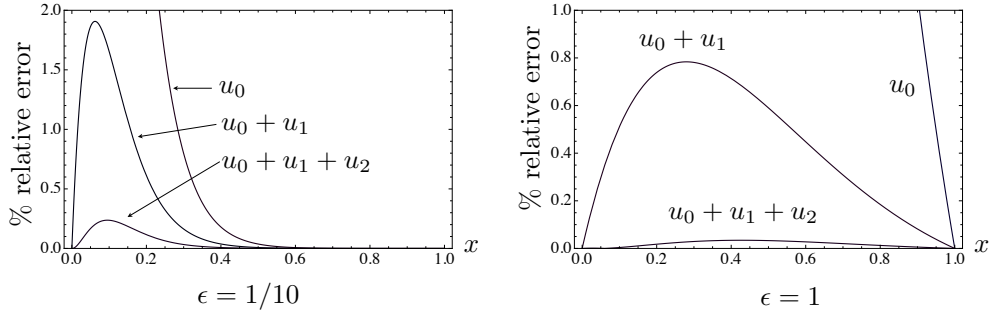}
    \caption{Comparison of relative percentage errors for $\epsilon=1/10$, $1$ with $\alpha=3/2, \beta=2$ for resummed approximations using only  $\tilde{u}_0(x)$, $\tilde{u}_0(x)+\tilde{u}_1(x)$ and $\tilde{u}_0(x)+\tilde{u}_1(x)+\tilde{u}_2(x)$.  Note the satisfaction of the boundary condition at both ends and the rapid apparent convergence.}
    \label{comp2}
       \end{center}
\end{figure}

At order ${\mathcal O}(\lambda)$ we then have
\begin{eqnarray}\label{pw5}
\qquad \epsilon  \tilde{u}_1''(x)+(x+\beta -1)\tilde{u}_1'(x)&=&0.
\end{eqnarray}
Applying the boundary condition, e.g., that  (with $\lambda=1$)
\begin{equation}\label{pw6}
\tilde{u}_0(0)+\tilde{u}_1(0)=\alpha, \qquad \tilde{u}_0(1)+\tilde{u}_1(1)=\beta ,
\end{equation}
we find
\begin{equation}\label{pw7}
\tilde{u}_1(x)=(\alpha -\beta +1) \left({\rm erf}\left(\frac{x+\beta -1}{\sqrt{2\epsilon } }\right)-{\rm erf}\left(\frac{\beta }{\sqrt{2\epsilon} }\right)\right) \slash \left({\rm erf}\left(\frac{\beta -1}{\sqrt{2\epsilon}}\right)-{\rm erf}\left(\frac{\beta }{\sqrt{2\epsilon }}\right)\right).
  \end{equation}
The numerical approximation of $\tilde{u}_0(x)+\tilde{u}_1(x)$ is better than any of the MAE or resummations so far (see figures \ref{comp1}, \ref{comp2}).  This is because the boundary conditions can here be satisfied at both ends by the leading two orders $\tilde{u}_0(x)$, $\tilde{u}_1(x)$.  Note that the exact forms (\ref{pw4}), (\ref{pw7}) must incorporate resummations of the $r$ series in algebraic powers of $\epsilon$, and this is here actually the reason why both boundary conditions have here been able to be satisfied.   That we can achieve this is due to the fortunate peculiar form of (\ref{pw5}) due to the truncation of the expansion in (\ref{pw4}).  This serendipity would not extend to general systems and the $r$ sum would not have been automatically resummed.

We can also solve the linear inhomogeneous equation at ${\mathcal O}(\lambda^2)$ for (\ref{node}) with boundary conditions $u_2(0)=u_2(1)=0$. Addition of $\tilde{u}_2(x)$ to  $\tilde{u}_0(x)+\tilde{u}_1(x)$ improves the numerical agreement even further (figure \ref{comp2}).

The form of $\tilde{u}_1(x)$ in relation to the double-transseries template requires an explanation.  If we substitute the asymptotic expansion
${\rm erf}(s)\sim 1-\exp \left(-s^2\right)/\sqrt{\pi } s$, $s \rightarrow +\infty$,
into $\tilde{u}_1(x)$ we find
\begin{equation}\label{pw10}
\tilde{u}_1(x)\sim  v(x)-v(1) \qquad 
v(x)=\frac{(\beta-1)(\alpha-\beta+1)\re^{-F(x)/\epsilon}}{(x+\beta-1)\left(1-(\beta-1)\re^{-F(1)/\epsilon}/\beta\right)},
\end{equation}
where $F(x)$ is as defined in the sections above.    At first sight, the form of $v(x)$ suggests that $\tilde{u}_1(x)$ might be interpreted as the  ${\mathcal O}(\epsilon^0)$ term of the resummed $p$-series for $n=1$,  (\ref{nresum}).  However a simple calculation shows this not to be the case.  In fact $\tilde{u}_1(x)$ does not satisfy the same boundary condition as the resummed (\ref{nresum}).  The reordering is more complicated with terms of ${\mathcal O}(\re^{-F(x)/\epsilon})$ appearing also in the higher order $\tilde{u}_k$.
Higher order $\tilde{u}_k$ satisfy inhomogeneous linear differential equations
\begin{equation}\label{pw15}
\epsilon  \tilde{u}_k''(x)+(x+\beta -1)\tilde{u}_k'(x)=f(\tilde{u}_0, \tilde{u}_0', \tilde{u}_1, \tilde{u}_1', \dots \tilde{u}_{k-1}),  \qquad \tilde{u}_k(0)=\tilde{u}_k(1)=0.
\end{equation}
The solutions of the homogeneous form equation will be of the same form for all $k$:
\begin{equation}\label{pw16}
A_k(\epsilon)v_1(x)+B_k(\epsilon)v_2(x)\re^{-F(x)/\epsilon}, \qquad v_l(x)\sim\sum_{r_{\rm min}}^\infty d_r^{(l)}(x)\epsilon^r, \ l=1,2.
\end{equation}
Terms of ${\mathcal O}(\exp(-jF(x)/\epsilon))$,  $0\le j\le k$ will arise from the inhomogeneous terms and so the $A_k(\epsilon)$ and $B_k(\epsilon)$  will contain terms in ${\mathcal O}(\re^{-jF(1)/\epsilon}), \ 0\le j\le k$ to satisfy the boundary conditions.  Hence the expansion (\ref{pw1}) can be seen as non-trivial reordered sum over $p$ and $r$, ordered by the index $n$.   For this particular system, due to the truncated nature of $\tilde{u}_0$, we are able to solve for the $\tilde{u}_k$ exactly and so effectively achieve this double resummation.  For situations where  $\tilde{u}_0$ does not truncate this will not be achievable in general.

Note that the transseries template is forced by the presence of boundary values, but the form of the lambda series is driven by the differential equations.  It is an open question as to what is the ordering of terms in the transseries that optimises the numerical agreement with the exact solution, when only a finite number of terms is taken in the approximation.

\section{Discussion}

We conclude with a brief comparison of our nonlinear results with those of Wasow (1956)
and a discussion of further work.  In our notation, Wasow considered (\ref{Wasoweq})
with $F_1(x,u,\epsilon)$ and $F_2(x,u,\epsilon)$ regular analytic with respect to $u$, $\epsilon$ and $C^{(2)}$ in $x$ in a region of $(x,u,\epsilon)$ space containing $u=u_0(x)$, $0\le x \le 1$ and $\epsilon=0$. 
Wasow's expansion takes the form
$$u(x,\epsilon)=\sum_{r=0}^{\infty}v_r(x,\epsilon)\mu^r, \qquad v_r(x,\epsilon)=\re^{-F(x)/\epsilon}w_r(x,\epsilon), \ r>0,$$
where $\mu=\alpha-v_0(0,\epsilon)=\alpha-\beta+1$, $v(1,\epsilon)=\beta$, the $w_r(x,\epsilon)$ are bounded and $F(x)$ is as defined above.  
He proved that the $v_r$ expansion is uniformly and absolutely convergent with respect to $\epsilon$, $\mu$ and $x$ for $0\le\epsilon\le\epsilon_1$, $|\alpha-\beta+1|\le \mu_1$, $0\le x\le 1.$
Wasow then reordered the terms to give
$$u(x,\epsilon)=\sum_{m=0}^\infty c_m(x,\alpha,\epsilon)\re^{-mF(x)/\epsilon}$$
as uniformly and absolutely convergent.  He stated that the coefficients $c_m(x,\alpha, \epsilon)$ are regular analytic and possess a power series in $\epsilon$. 
He omited the full proof as it was, in his own words, ``somewhat detailed".

Wasow's results are existence proofs.  He took the resummed route directly, his ``reordering" is related to the transseries.  Here we have taken the reverse route of solving for the transseries first and then resumming as we have been motivated by an exponential asymptotics approach.  Our explicit calculation reveals the intricate and subtle nature of the form of Wasow's $c_m$ or $w_r$ coefficients and reveals an interesting realignment of exponential scales that are required to satisfy the boundary values. The template has been dictated by the boundary values and is independent of the equation and so is expected to be more generally valid.  Wasow's statements that the coefficients ``possess an asymptotic expansion in powers of $\epsilon$" is strictly incomplete, since it ignores the presence of series beyond all orders, and in inverse powers of the asymptotic parameter.   Given all this, Wasow's work suggests that our resummations may be convergent, at least for small values of $\mu=\alpha-\beta+1$.

The preliminary work of this paper has opened up a multitude of potential follow-up problems in exponential asymptotics. Further work could include:  rigorous justification of the exponential approach, with a proof (or otherwise) of convergence; derivation of proper error bounds for the reordered summations; extensions to higher order BVPs; adjustment of the template to include interior layers and shocks; a full hyperasymptotic treatment based on the fundamental transseries template;  extensions of the transeries template to PDE BVPs.  Some of these will be discussed elsewhere.

\begin{acknowledgements}

The author acknowledges the hospitality of the Pacific Institute for the Mathematical Sciences and the Department of Mathematics at UBC, Vancouver, where some of this work was undertaken, and AB Olde Daalhuis for very helpful discussions.

\end{acknowledgements}

\appendix{Multiple-scales exponential asymptotics}

Motivated by the form of (\ref{psum}) we seek, {\it a priori}, a multiple-scales solution of (\ref{1}) that corresponds to the summation of all the $p$-exponentials in (\ref{template}).  We use the ansatz:
\begin{equation}\label{Wansatz}
u(x;\epsilon)\sim \sum_{r=0}^{\infty}W_r(x,X)\epsilon^r+\re^{-F(x)/\epsilon}\sum_{r=0}^{\infty}V_n(x,X)\epsilon^r,
\end{equation}
where the scaled ``variable" is the constant $X=F(1)/\epsilon$.

We substitute(\ref{Wansatz}) into (\ref{1}) and balance at ${\mathcal O}(\epsilon^r)$ and ${\mathcal O}(\re^{-F(x)/\epsilon}\epsilon^r)$, ignoring the  $\epsilon$-dependence in the $X$ terms.  Since $X$ is a constant, it does not actually generate a derivative in $\partial/\partial X$.  We thus obtain recurrence relations that are identical to those we obtained when we substituted (\ref{latta}) into (\ref{1}). Hence the recurrence relations for W and V are, for $r\ge 0$,
\begin{equation}\label{Wc}
W_r(x,X)=(c_r-W_{r-1}'(x,X))/(2x+1), \ \ V'_r(x,X)=(V_{r-1}''(x,X))/(2x+1),
\end{equation}
with $W_{-1}(x,X)=$ $V_{-1}(x,X)$ $=0$ and prime denoting $x$-differentiation.

The constants $c_r$  are determined from a modified set of boundary conditions.  The terms in the expansion in (\ref{bc0}), (\ref{bc1}) and (\ref{pbcs0}) have been balanced at orders of ${\mathcal O}(\re^{-F(x)/\epsilon}\epsilon^r)$.  Starting from (\ref{Wansatz}) due to the apparent similarity of (\ref{latta}) and the second row of (\ref{Wansatz}), we find that we go around in a circle.  The second sum in (\ref{Wansatz}) cannot satisfy the boundary condition at $x=1$ exactly, without including further series, for the reasons outlined above.  Thus the boundary conditions have to be modified.  Instead of (\ref{bc0}-\ref{bc1}) we have:
\begin{eqnarray}\label{Wbc0}
W_r(0,X)+V_r(0,X)=\delta_{r0}\alpha, \ \ W_r(1,X)+\re^{-X}V_r(1,X)=\delta_{r0}\beta. 
\end{eqnarray}
Note that the exponential in (\ref{Wbc0}) is actually $\re^{-F(1)/\epsilon}$, precisely the order of neglected terms in (\ref{bc1}) which lead to the ansatz (\ref{template}). That we can include this exponential term now is because we treat $X$ as varying on a different scale to $\epsilon$.  

From (\ref{Wbc0}) the $c_r$ in (\ref{Wc}) can now be found.  A short calculation gives
\begin{eqnarray}\label{Wr} \nonumber
W_r(x,X)&=&\frac{(W_{r-1}'(1,X)-W_{r-1}'(x,X))-3\re^{-X}(W'_{r-1}(0,X)-W'_{r-1}(x,X))}{(1-3\re^{-X})(2x+1)}, \\ \label{Vr}
V_r(x,X)&=&(W_{r-1}'(0)-W_{r-1}'(0))/(1-3\re^{-X}).
\end{eqnarray}
with initial terms $W_0(x,X)=3(\beta-\alpha\re^{-X})/\left\{(1-3\re^{-X})(2x+1)\right\}$ and 
$V_0(x,X)=(\alpha-3\beta)/(1-3\re^{-X})$. Inserting the value of $\epsilon X=F(1)=2$, to leading order in $\epsilon$ we recover (\ref{psum}).  Note that by neglecting the terms in $\re^{-X}$ these relations reduce to those of (\ref{aterm}-\ref{bterm}), as they should.  


\label{lastpage}

 \end{document}